# SOME RESULTS ON TWO-SIDED LIL BEHAVIOR

By Uwe Einmahl[1] and Deli Li[2]

*Vrije Universiteit Brussel and Lakehead University*

Let $\{X, X_n; n \geq 1\}$ be a sequence of i.i.d. mean-zero random variables, and let $S_n = \sum_{i=1}^{n} X_i, n \geq 1$. We establish necessary and sufficient conditions for having with probability 1, $0 < \limsup_{n\to\infty} |S_n|/\sqrt{nh(n)} < \infty$, where $h$ is from a suitable subclass of the positive, nondecreasing slowly varying functions. Specializing our result to $h(n) = (\log\log n)^p$, where $p > 1$ and to $h(n) = (\log n)^r$, $r > 0$, we obtain analogues of the Hartman–Wintner LIL in the infinite variance case. Our proof is based on a general result dealing with LIL behavior of the normalized sums $\{S_n/c_n; n \geq 1\}$, where $c_n$ is a sufficiently regular normalizing sequence.

**1. Introduction.** Let $\{X, X_n; n \geq 1\}$ be a sequence of real-valued independent and identically distributed (i.i.d.) random variables, and let $S_n = \sum_{i=1}^{n} X_i, n \geq 1$. Define $Lx = \log_e \max\{e, x\}$ and $LLx = L(Lx)$ for $x \in \mathbb{R}$. The classical Hartman–Wintner law of the iterated logarithm (LIL) states that

$$\limsup_{n\to\infty} \pm S_n/(2nLLn)^{1/2} = \sigma \qquad \text{a.s.} \tag{1.1}$$

if and only if

$$\mathbb{E}X = 0 \quad \text{and} \quad \sigma^2 = \mathbb{E}X^2 < \infty. \tag{1.2}$$

Moreover, if (1.2) holds, then

$$C(\{S_n/\sqrt{2nLLn}; n \geq 1\}) = [-\sigma, \sigma] \qquad \text{a.s.,} \tag{1.3}$$

where $C(\{x_n; n \geq 1\})$ stands for the cluster set (i.e., the set of limit points) of the sequence $\{x_n; n \geq 1\}$. See [8] for the "if" part and [24] for the converse.

Received July 2004; revised October 2004.
[1]Supported in part by an FWO grant.
[2]Supported by a grant from the Natural Sciences and Engineering Research Council of Canada.
*AMS 2000 subject classifications.* 60F15, 60G50.
*Key words and phrases.* Hartman–Wintner LIL, law of the iterated logarithm, superslow variation, two-sided LIL behavior, sums of i.i.d. random variables, cluster sets.

The conclusion (1.3) is due to Strassen [23]. Actually, in this fundamental paper, Strassen [23] obtained a functional LIL as well as invariance principles which are in many respects at the origin of the study of LIL in a vector-valued setting. For very efficient and self-contained proofs of the Hartman–Wintner LIL which do not use the Kolmogorov LIL [14] see, for example, [2] or [7].

It is natural then to ask whether one can find analogous results for variables with infinite variance. This of course requires different normalizing sequences and also sometimes different centering sequences. In the case where $\{X, X_n; n \geq 1\}$ is a sequence of symmetric real-valued i.i.d. random variables, Feller [6] (see [10, 19] and [3] for some clarification) studied the problem of determining when it is possible to find a positive regular monotone sequence $\{a_n; n \geq 1\}$ such that

$$\limsup_{n\to\infty} |S_n|/a_n = 1 \qquad \text{a.s.} \tag{1.4}$$

In this case, one speaks of *two-sided* LIL behavior.

Of course one can also address the corresponding *one-sided* LIL behavior problem with centering $\{\delta_n\}$, that is, when one has for a suitable (regular) sequence $b_n$

$$0 < \limsup_{n\to\infty} (S_n - \delta_n)/b_n < \infty \qquad \text{a.s.} \tag{1.5}$$

For some basic work in this direction refer to [11, 12, 13] in the finite expectation case where $\delta_n = n\mathbb{E}X$ and for more general results see also [22].

Kuelbs and Zinn [17] showed that the techniques of Klass [11] are also extremely useful for the LIL problem in Banach space, and this was further elaborated by Kuelbs [16] and Einmahl [4]. The main purpose of the present paper is to address some still open questions in connection with two-sided LIL behavior for real-valued random variables with finite expectation.

To cite the relevant work in this direction let us first recall some definitions of Klass [11]. As above, let $X : \Omega \to \mathbb{R}$ be a random variable and assume that $0 < \mathbb{E}|X| < \infty$. Set

$$H(t) := \mathbb{E}X^2 I\{|X| \leq t\} \quad \text{and} \quad M(t) := \mathbb{E}|X| I\{|X| > t\}, \qquad t \geq 0.$$

Then it is easy to see that the function $G(t) := t^2/(H(t) + tM(t))$, $t > 0$, is continuous and increasing with an inverse function $K(x)$, $x > 0$. Moreover, one has for this function $K$ that as $x \nearrow \infty$

$$K(x)/\sqrt{x} \nearrow (\mathbb{E}X^2)^{1/2} \in ]0, \infty] \tag{1.6}$$

and

$$K(x)/x \searrow 0. \tag{1.7}$$

Set $\gamma_n = \sqrt{2}K(n/LLn)LLn$, $n \geq 1$. Klass [11, 12] established a one-sided LIL result with respect to this sequence which also implies the following two-sided LIL result if $\mathbb{E}X = 0$:

$$\limsup_{n\to\infty} |S_n|/\gamma_n = 1 \tag{1.8}$$

if and only if

$$\sum_{n=1}^{\infty} \mathbb{P}\{|X| \geq \gamma_n\} < \infty. \tag{1.9}$$

(Actually, Klass [11, 12] proved that the limiting constant is $\in [2^{-1/2}, 3 \cdot 2^{-3/2}]$ and showed later in [13] that this is optimal in the one-sided case. For the calculation of the limiting constant in the two-sided case, see [4] and also Section 3 of the present paper.)

We thus see that if condition (1.9) is satisfied, one obtains an LIL result which extends the classical Hartman–Wintner LIL. [Note that if $0 < \sigma^2 = \mathbb{E}X^2 < \infty$, we have that $K(n/LLn)LLn \sim \sigma\sqrt{nLLn}$ and condition (1.9) is trivially satisfied so that the Hartman–Wintner LIL is a special case of (1.8).]

Moreover, Klass [11] (see his Theorem 4.1) has shown that if $\mathbb{E}X = 0$ and $c_n \geq \sqrt{9/8}\gamma_n$ is a sequence so that $c_n/n^{1/2}$ is increasing, one has

$$\limsup_{n\to\infty} |S_n|/c_n \leq 1 \qquad \text{a.s.} \tag{1.10}$$

if and only if

$$\sum_{n=1}^{\infty} \mathbb{P}\{|X| \geq c_n\} < \infty. \tag{1.11}$$

This implies for sequences $c_n$ satisfying $c_n/\gamma_n \to \infty$

$$\limsup_{n\to\infty} |S_n|/c_n = 0 \text{ or } \infty \qquad \text{a.s.} \tag{1.12}$$

according as

$$\sum_{n=1}^{\infty} \mathbb{P}\{|X| \geq c_n\} < \infty \qquad \text{or } = \infty. \tag{1.13}$$

We thus see that if one considers "big" sequences as above, one can only obtain stability results, but no longer LIL behavior.

Here we shall investigate whether there are still LIL type results if condition (1.9) is not satisfied and, moreover, whether one can find "nicer" norming sequences than $\{\gamma_n\}$. This sequence is very appealing in that it is defined in a universal way depending on the distribution of $X$ only, but if one looks at concrete examples it can be quite difficult to determine $\{\gamma_n\}$. Another problem is that in certain situations the sequence $\gamma_n$ can be too small. An example which was discussed by Feller [6] and Pruitt [22] is a symmetric random variable $X$ with Lebesgue density $f_X(x) = |x|^{-3}$, $|x| \geq 1$. In this case it is easy to calculate $\gamma_n$, but assumption (1.11) is not satisfied so that the LIL of Klass does not apply, nor do the LIL results of Feller [6] and Pruitt [22]. It seems to be still an open problem whether, in this particular case, there exists a "nice" normalizing sequence $a_n$ so that $\limsup_{n\to\infty} |S_n|/a_n = 1$ a.s.

We first address the following modified form of the LIL behavior problem.

Problem 1. Given a sequence, $a_n = \sqrt{nh(n)}$, where $h$ is a slowly varying nondecreasing function, we ask: When do we have with probability 1, $0 < \limsup_{n\to\infty} |S_n|/a_n < \infty$?

One possibility would be to look for conditions implying $\gamma_n \approx a_n$, but as we are dealing with almost sure convergence one has many more possibilities for finding normalizing sequences than in the weak convergence case. Under an additional assumption on $h$ we will establish a necessary and sufficient condition for LIL behavior with respect to the given sequence $a_n$. Using this result we can also find a normalizing sequence of this type for the Feller–Pruitt example (see Section 5 below).

At first sight our result might look quite different from the Klass LIL, but it will turn out that our conditions imply

$$0 < \liminf_{n\to\infty} a_n/\gamma_n < \infty \tag{1.14}$$

which shows that we are in the range between the LIL result (1.8) and the stability result (1.12). It is natural then to pose a second related question, namely

Problem 2. Consider a nondecreasing sequence $c_n$ satisfying $0 < \liminf_{n\to\infty} c_n/\gamma_n < \infty$. When do we have with probability 1, $0 < \limsup_{n\to\infty} |S_n|/c_n < \infty$? If this is the case, what is the cluster set $C(\{S_n/c_n; n \geq 1\})$?

From Corollary 10 of [20] in combination with (3.5) below it follows that, under a mild regularity assumption on the sequence $\{c_n\}$, the above lim sup is equal to a certain parameter $\alpha_0$. We shall additionally show that the corresponding cluster set $C(\{S_n/c_n; n \geq 1\})$ always coincides with the interval $[-\alpha_0, \alpha_0]$ (see Theorem 3 below). It is then clear that we have LIL behavior with respect to the normalizing sequence $c_n$ if and only if $0 < \alpha_0 < \infty$. Thus, in principle, this solves Problem 2. There is still a difficulty, namely, the determination of this parameter. For that reason, we shall also show that under assumption (1.14) one can define this parameter differently, which makes the calculation of $\alpha_0$ feasible in many cases of interest (see Theorem 4). This way we can immediately reobtain the two-sided version of the Klass LIL (1.8) and we get a whole class of new LIL type results as indicated in Problem 1. (For a survey of some other work on Problem 2 refer to Sections 7.3 and 7.5 of [21].)

The plan of the paper is as follows. Our main results regarding Problem 1, Theorems 1 and 2, and their corollaries as well as Theorems 3 and 4 are presented in Section 2. In Section 3 we prove the two latter theorems and in Section 4 we show how one can infer Theorems 1 and 2 from them. After giving a few examples and some further comments in Section 5, we finally determine the desired "nice" normalizing sequence for the Feller–Pruitt example.

**2. Statement of main results.** Before we can formulate our results, we need some extra notation. Let $\mathcal{H}$ be the set of all continuous, nondecreasing functions $h:[0,\infty[\to]0,\infty[$, which are slowly varying at infinity. By monotonicity the slow variation of $h$ is equivalent to $\lim_{t\to\infty} h(et)/h(t)=1$. Very often one can even show that $\lim_{t\to\infty} h(tf(t))/h(t)=1$, where $f$ is an increasing function such that $\lim_{t\to\infty} f(t)=\infty$. For instance, if $h(t)=LLt$, $t\geq 0$, this is the case for $f(t)=t$. In the literature this is also called super-slow variation (refer to pages 186–188 in [1] for more information and background on this notion).

For our purposes the functions $f_\tau(t):=\exp((Lt)^\tau)$, $0\leq\tau\leq 1$, will be most important. Clearly if $\lim_{t\to\infty} h(tf(t))/h(t)=1$ holds for $f=f_\tau$, where $\tau>0$ this also holds for $f=f_{\tau'}$, $0\leq\tau'\leq\tau$. Thus, the bigger we can choose the parameter $\tau$, the slower is the variation of the given function $h$. (Also note that this condition with $\tau=0$ is equivalent with slow variation.)

Given $0\leq q<1$, let $\mathcal{H}_q\subset\mathcal{H}$ be the class of all functions so that

$$\lim_{t\to\infty} h(tf_\tau(t))/h(t)=1, \qquad 0<\tau<1-q,$$

and set $\mathcal{H}_1=\mathcal{H}$. We consider $q$ as a measure for how slow the variation is. So functions in $\mathcal{H}_0$ are the "slowest" and it will turn out that this class is particularly interesting for LIL type results (see Theorem 2 below). Examples for functions in $\mathcal{H}_0$ are $h(t)=(Lt)^r$, $r\geq 0$, and $h(t)=(LLt)^p$, $p\geq 0$.

The following Theorem 1 gives LIL type results if $\lambda>0$ and stability results if $\lambda=0$ with respect to a large class of normalizing sequences, without assuming that $\mathbb{E}X^2<\infty$.

THEOREM 1. *Let $X,X_1,X_2,\ldots$ be i.i.d. random variables, and let $S_n=\sum_{i=1}^n X_i, n\geq 1$. Given a function $h\in\mathcal{H}_q$ where $0\leq q\leq 1$, set $\Psi(x)=\sqrt{xh(x)}$ and $a_n=\Psi(n)$, $n\geq 1$. If there exists a constant $0\leq\lambda<\infty$ such that*

$$\mathbb{E}X=0, \qquad \mathbb{E}\Psi^{-1}(|X|)<\infty, \qquad \limsup_{x\to\infty}\frac{\Psi^{-1}(xLLx)}{x^2LLx}H(x)=\frac{\lambda^2}{2}, \tag{2.1}$$

*then we have*

$$(1-q)^{1/2}\lambda \le \limsup_{n\to\infty} |S_n|/a_n \le \lambda \qquad a.s. \tag{2.2}$$

*Conversely, if $q<1$, then the relation*

$$\limsup_{n\to\infty} \frac{|S_n|}{a_n} < \infty \qquad a.s. \tag{2.3}$$

*implies that (2.1) holds for some $\lambda < \infty$.*

*Moreover, the* lim sup *in (2.3) is positive if and only if (2.1) holds for some $\lambda > 0$.*

Note the lim sup in condition (2.1). If this is actually a limit or if the corresponding lim inf is positive, one can show that $a_n \approx \gamma_n$ and one could obtain (2.3) from the Klass LIL (with less tight bounds on the limiting constant). This is no longer possible if the lim inf is equal to 0, which clearly indicates that we can obtain LIL type results in many situations where the Klass LIL does not apply. The reader will notice that we have taken advantage of this additional possibility for proving such results when choosing $a_n$ in the Feller–Pruitt example (see Section 5).

For slowly varying functions $h \in \mathcal{H}_0$ we obtain a complete analogue of the Hartman–Wintner LIL.

THEOREM 2. *Assume that $h \in \mathcal{H}_0$ and let $\Psi$ and $\{a_n\}$ be as in Theorem 1. For any constant $0 \le \lambda < \infty$ we have:*

$$\limsup_{n\to\infty} \pm S_n/a_n = \lambda \qquad a.s. \tag{2.4}$$

*and*

$$C(\{S_n/a_n; n \ge 1\}) = [-\lambda, \lambda] \qquad a.s. \tag{2.5}$$

*if and only if condition (2.1) holds.*

We shall illustrate Theorem 2 by considering the following two special cases:

*Case* 1. Take $h(x) = 2(LLx)^p$ where $p \ge 1$. Then it is easy to check that

$$\lim_{x\to\infty} \frac{\Psi^{-1}(x)}{x^2/(2(LLx)^p)} = 1.$$

It follows that

$$\lim_{x\to\infty} \frac{\Psi^{-1}(xLLx)/(x^2LLx)}{1/(2(LLx)^{p-1})} = 1.$$

*Case* 2. Choose $h(x) = 2(Lx)^r$ where $r > 0$. One easily sees that

$$\lim_{x\to\infty} \frac{\Psi^{-1}(x)}{x^2/(Lx)^r} = 2^{-(r+1)}$$

and

$$\lim_{x\to\infty} \frac{\Psi^{-1}(xLLx)/(x^2LLx)}{LLx/(Lx)^r} = 2^{-(r+1)}.$$

Thus Theorem 2 implies the following two results.

Corollary 1. *Let $p \geq 1$. For any constant $0 \leq \lambda < \infty$ we have:*

$$\limsup_{n\to\infty} \frac{\pm S_n}{\sqrt{2n(LLn)^p}} = \lambda \qquad a.s.$$

*if and only if*

$$\text{(2.6)}\quad \mathbb{E}X = 0, \qquad \mathbb{E}X^2/(LL|X|)^p < \infty, \qquad \limsup_{x\to\infty} (LLx)^{1-p} H(x) = \lambda^2.$$

Remark 1. If $p = 1$, then condition (2.6) is equivalent to

$$\mathbb{E}X = 0 \quad \text{and} \quad \mathbb{E}X^2 = \lambda^2.$$

We see that the classical Hartman–Wintner LIL is a special case of Corollary 1.

Corollary 2. *Let $r > 0$. For any constant $0 \leq \lambda < \infty$ we have:*

$$\limsup_{n\to\infty} \frac{\pm S_n}{\sqrt{2n(Ln)^r}} = \lambda \qquad a.s.$$

*if and only if*

$$\text{(2.7)}\quad \mathbb{E}X = 0, \qquad \mathbb{E}X^2/(L|X|)^r < \infty, \qquad \limsup_{x\to\infty} \frac{LLx}{(Lx)^r} H(x) = 2^r\lambda^2.$$

For a further corollary to Theorem 1 (where $0 < q < 1$) refer to Section 5. If condition (2.1) in Theorem 1 is satisfied with $\lambda = 0$ we obtain the following stability result.

Corollary 3. *Let $h \in \mathcal{H}$ and let $\Psi$ and $\{a_n\}$ be as in Theorem 1. If*

$$\text{(2.8)}\quad \mathbb{E}X = 0, \qquad \mathbb{E}\Psi^{-1}(|X|) < \infty, \qquad \lim_{x\to\infty} \frac{\Psi^{-1}(xLLx)}{x^2LLx} H(x) = 0,$$

*then*

$$\text{(2.9)}\qquad \lim_{n\to\infty} S_n/a_n = 0 \qquad a.s.$$

*Moreover, if $h \in \mathcal{H}_q$ for some $q < 1$, then condition (2.8) is necessary and sufficient for (2.9) to hold.*

REMARK 2. We note that after some work (2.9) also follows from (1.12) (see Remark 5 in Section 4). The necessity of condition (2.8) is a new result as far as we know.

We first look at Problem 2 for sequences $c_n$ satisfying the following two conditions:

$$c_n/\sqrt{n} \nearrow \infty \tag{2.10}$$

and

$$\forall \varepsilon > 0 \; \exists m_\varepsilon \geq 1 : c_n/c_m \leq (1+\varepsilon)(n/m), \qquad m_\varepsilon \leq m < n. \tag{2.11}$$

Note that condition (2.11) is satisfied if $c_n/n$ is nonincreasing (e.g., if $c_n = \gamma_n$) or if $c_n = c(n)$, where $c : [0, \infty) \to [0, \infty)$ is regularly varying at infinity with exponent $\gamma < 1$. (This includes all the sequences $\{a_n\}$ considered in Problem 1.)

THEOREM 3. *Let $X, X_1, X_2, \ldots$ be i.i.d. mean-zero random variables. Assume that*

$$\sum_{n=1}^{\infty} \mathbb{P}\{|X| \geq c_n\} < \infty, \tag{2.12}$$

*where $c_n$ is a sequence of positive real numbers satisfying conditions* (2.10) *and* (2.11). *Set*

$$\alpha_0 = \sup\left\{\alpha \geq 0 : \sum_{n=1}^{\infty} n^{-1} \exp\left(-\frac{\alpha^2 c_n^2}{2n\sigma_n^2}\right) = \infty\right\},$$

*where $\sigma_n^2 = H(\delta c_n)$ and $\delta > 0$.*

*Then we have with probability* 1,

$$C(\{S_n/c_n; n \geq 1\}) = [-\alpha_0, \alpha_0] \tag{2.13}$$

*and*

$$\limsup_{n\to\infty} |S_n|/c_n = \alpha_0. \tag{2.14}$$

REMARK 3. As mentioned above, (2.14) also follows from Corollary 10 of [20], where the parameter $\alpha_0$ has been defined slightly differently. It is easy to see that our definition is consistent with his definition. Also note that $\alpha_0$ can be infinite. [Choose, e.g., $c_n = n^{1/2}(LLn)^{1/4}$.]

THEOREM 4. *Let $X$ and $c_n$ be as in Theorem* 3. *Further assume that $a := \liminf_{n\to\infty} c_n/\gamma_n > 0$. Then we can choose $\sigma_n^2$ in the definition of $\alpha_0$ equal to $H(d_n)$, where $d_n \leq c_n$ can be any sequence satisfying*

$$\log(c_n/d_n)/LLn \to 0 \qquad \text{as } n \to \infty. \tag{2.15}$$

*Moreover, we have in this case $\alpha_0 \leq 1/a < \infty$.*

Remark 4. Note that Theorem 4 also gives the upper bound part of the LIL result (1.8) (just set $c_n = \gamma_n$). In general, this result will be very helpful for finding upper bounds for $\alpha_0$ as it allows us to replace $\delta c_n$ by a "small" $d_n$. If one wants to find a lower bound for $\alpha_0$ one normally should choose $d_n = c_n$, and Theorem 3 will be sufficient. So it is not too surprising that the lower bound part of (1.8) already follows from Theorem 3 (see end of Section 3).

**3. Proofs of Theorems 3 and 4.** Throughout the whole section we assume that $\{c_n\}$ is a sequence of positive real numbers satisfying conditions (2.10) and (2.11). Moreover, $X, X_1, X_2, \ldots$ will always be a sequence of i.i.d. mean-zero random variables satisfying

$$\sum_{n=1}^{\infty} \mathbb{P}\{|X| \geq c_n\} < \infty.$$

In the first lemma we collect some more or less known facts which we need for the proof of Theorem 3.

Lemma 1. *We have*

$$\sum_{n=1}^{\infty} \mathbb{E}|X|^3 I\{|X| \leq c_n\}/c_n^3 < \infty, \tag{3.1}$$

$$\sum_{n=1}^{\infty} \mathbb{P}\{|X| > \varepsilon c_n\} < \infty \qquad \forall\, \varepsilon > 0, \tag{3.2}$$

$$H(c_n) = \mathbb{E}X^2 I\{|X| \leq c_n\} = o(c_n^2/n) \qquad \textit{as } n \to \infty, \tag{3.3}$$

$$M(c_n) = \mathbb{E}|X| I\{|X| > c_n\} = o(c_n/n) \qquad \textit{as } n \to \infty, \tag{3.4}$$

$$\mathbb{E}|S_n| = o(c_n) \qquad \textit{as } n \to \infty. \tag{3.5}$$

Proof. For the first fact refer, for instance, to Lemma 1 of [3]. We only need to prove (3.2) if $\varepsilon < 1$. In this case it directly follows from (3.1) via the inequality

$$\mathbb{P}\{\varepsilon c_n \leq |X| < c_n\} \leq \varepsilon^{-3} \mathbb{E}|X|^3 I\{|X| \leq c_n\}/c_n^3.$$

To prove (3.3) we first note that $\sum_{n=1}^{\infty} \mathbb{P}\{|X| \geq c_n\} < \infty$ is equivalent to $\sum_{j=1}^{\infty} j p_j < \infty$, where $p_j = \mathbb{P}\{c_{j-1} < |X| \leq c_j\}$, $j \geq 1$ (with $c_0 = 0$).

Then we readily obtain for any $j_0 \geq 1$ and $n \geq j_0 + 1$,

$$nH(c_n)/c_n^2 = nH(c_{j_0})/c_n^2 + n \sum_{j=j_0+1}^{n} \{H(c_j) - H(c_{j-1})\}/c_n^2$$

$$\leq nH(c_{j_0})/c_n^2 + \sum_{j=j_0+1}^{n} p_j n(c_j/c_n)^2$$

$$\leq nH(c_{j_0})/c_n^2 + \sum_{j=j_0+1}^{\infty} jp_j.$$

Choosing $j_0$ so large that $\sum_{j=j_0+1}^{\infty} jp_j < \varepsilon$, we see that

$$\limsup_{n\to\infty} nH(c_n)/c_n^2 \leq \varepsilon, \qquad \varepsilon > 0,$$

which proves (3.3).

To see (3.4) simply note that on account of (2.11) there exists a constant $K \geq 1$ so that

$$n\mathbb{E}|X|I\{|X| > c_n\}/c_n \leq n \sum_{j=n+1}^{\infty} c_j p_j/c_n \leq K \sum_{j=n+1}^{\infty} jp_j,$$

which goes to zero as $n \to \infty$.

If $X$ has a symmetric distribution we have

$$\mathbb{E}|S_n| \leq (nH(c_n))^{1/2} + n\mathbb{E}|X|I\{|X| > c_n\},$$

and fact (3.5) follows in this case by combining the two previous facts. Using a standard symmetrization argument, we obtain (3.5) for nonsymmetric random variables as well. □

We now determine the cluster set $C(\{S_n/c_n; n \geq 1\}) =: A$, where we use Theorem 3 of [9]. (It is easily seen that $c_n$ satisfies the conditions of this result.) Since $S_n/c_n \xrightarrow{\mathbb{P}} 0$ [see (3.5)], it follows from Kesten's result (see also [15]) that

$$x \in C(\{S_n/c_n; n \geq 1\}) \iff \sum_{n=1}^{\infty} \frac{1}{n}\mathbb{P}\{|S_n/c_n - x| < \varepsilon\} = \infty \quad \forall \varepsilon > 0. \tag{3.6}$$

Using this equivalence, one can further prove

LEMMA 2. *We have*

$$x \in A \iff \sum_{n=1}^{\infty} \frac{1}{n}\mathbb{P}\{|S_{n,n}/c_n - x| < \varepsilon\} = \infty \qquad \forall \varepsilon > 0, \tag{3.7}$$

*where* $S_{n,n} = \sum_{i=1}^{n}\{X_{n,i} - \mathbb{E}X_{n,i}\}$, $X_{n,i} = X_i I\{|X_i| \leq d_n\}$ *and* $d_n = \delta c_n$, *with* $\delta > 0$.

PROOF. In view of (3.6) it is enough to show that

$$\sum_{n=1}^{\infty} n^{-1}\mathbb{P}\{|S_n - S_{n,n}| \geq \varepsilon c_n\} < \infty \qquad \forall\, \varepsilon > 0. \tag{3.8}$$

Recalling (3.2) and (3.4), we have as $n \to \infty$,

$$\left|\sum_{i=1}^{n} \mathbb{E}X_{n,i}\right| \leq n\mathbb{E}|X|I\{|X| > \delta c_n\} = o(c_n)$$

and we can infer that for large $n$

$$\mathbb{P}\{|S_n - S_{n,n}| \geq \varepsilon c_n\} \leq \mathbb{P}\left\{S_n \neq \sum_{i=1}^{n} X_{n,i}\right\}$$

which is less than or equal to

$$n\mathbb{P}\{|X| > \delta c_n\}$$

and we readily obtain (3.7) from (3.2). □

From (3.5) we obtain that $0 \in A$ and we can focus on the nonzero elements of $A$.

LEMMA 3. *Let $x \neq 0$. Then we have*

$$x \in A \quad \Longleftrightarrow \quad \sum_{n=1}^{\infty} \frac{1}{n}\mathbb{P}\{|\sigma_n\sqrt{n}Z/c_n - x| < \varepsilon\} = \infty \qquad \forall\, \varepsilon > 0, \tag{3.9}$$

*where $Z$ is a standard normal variable and $\sigma_n^2 = H(\delta c_n)$, with $\delta > 0$.*

PROOF. Using a well-known nonuniform Berry–Esseen type inequality (see, e.g., Theorem 5.17 on page 168 in [21]), it follows that for $0 < \varepsilon < |x|/2$,

$$\begin{aligned} &|\mathbb{P}\{|S_{n,n}/c_n - x| < \varepsilon\} - \mathbb{P}\{|\tilde{\sigma}_n\sqrt{n}Z/c_n - x| < \varepsilon\}| \\ &\qquad \leq 16C|x|^{-3}n\mathbb{E}|X_{n,1} - \mathbb{E}X_{n,1}|^3/c_n^3 \leq 128C|x|^{-3}n\mathbb{E}|X|^3 I\{|X| \leq c_n\}/c_n^3, \end{aligned}$$

where in the last step we have used the $c_r$-inequality. $C$ is an absolute constant and $\tilde{\sigma}_n^2 = \mathrm{Var}(X_{n,1})$. Recalling (3.1) we see that $x \in A$ is equivalent to

$$\sum_{n=1}^{\infty} \frac{1}{n}\mathbb{P}\{|\tilde{\sigma}_n\sqrt{n}Z/c_n - x| < \varepsilon\} = \infty \qquad \forall\, \varepsilon > 0. \tag{3.10}$$

Let $\delta_n^2 = \sigma_n^2 - \tilde{\sigma}_n^2 = (\mathbb{E}X_{n,1})^2$. By the dominated convergence theorem we have $\delta_n \to 0$ as $n \to \infty$ and recalling that $\sigma_n^2 \nearrow \mathbb{E}X^2 > 0$ we see that

$$\sigma_n^2/\tilde{\sigma}_n^2 \to 1 \qquad \text{as } n \to \infty, \tag{3.11}$$

from which we readily obtain that the series condition in (3.7) is equivalent to (3.9) and the lemma has been proven. □

Using the trivial inequality $\mathbb{P}\{Z > t+s\} \le \mathbb{P}\{Z > t\}/2$, $s > 1/t$, $t > 0$, we can further simplify the lemma about clustering as follows.

LEMMA 4. *We have,*

$$x \in A \quad \Longleftrightarrow \quad \sum_{n=1}^{\infty} \frac{1}{n} \exp\left(-\frac{(|x|-\varepsilon)_+^2 c_n^2}{2n\sigma_n^2}\right) = \infty \qquad \forall\, \varepsilon > 0, \tag{3.12}$$

*where $\sigma_n^2$ is defined as in Lemma 3.*

PROOF. If $x = 0$, the equivalence is trivial. If $x > 0$, we have in view of Lemma 3 that $x \in A$ is equivalent to

$$\sum_{n=1}^{\infty} \frac{1}{n} \mathbb{P}\{x - \varepsilon < \sigma_n \sqrt{n} Z/c_n < x + \varepsilon\} = \infty \qquad \forall\, \varepsilon > 0. \tag{3.13}$$

This trivially implies that

$$\sum_{n=1}^{\infty} \frac{1}{n} \mathbb{P}\{x - \varepsilon < \sigma_n \sqrt{n} Z/c_n\} = \infty \qquad \forall\, \varepsilon > 0, \tag{3.14}$$

which in turn by standard estimates of the tail probabilities of the normal distribution is equivalent to the series condition in (3.12). It remains to show that (3.14) implies (3.13). To that end we note that if $\varepsilon < x/2$,

$$\mathbb{P}\{\sigma_n \sqrt{n} Z/c_n > x+\varepsilon\} = \mathbb{P}\left\{Z > \frac{(x-\varepsilon)c_n}{\sqrt{n}\sigma_n} + \frac{2\varepsilon c_n}{\sqrt{n}\sigma_n}\right\} \le \frac{1}{2}\mathbb{P}\left\{Z > \frac{(x-\varepsilon)c_n}{\sqrt{n}\sigma_n}\right\},$$

provided that

$$2\varepsilon c_n/(\sqrt{n}\sigma_n) \ge \sqrt{n}\sigma_n/\{(x-\varepsilon)c_n\}.$$

Relation (3.3) implies $\sigma_n^2 = H(\delta c_n) = o(c_n^2/n)$ and it follows that the above condition is satisfied for large $n$. We thus have in this case,

$$2\mathbb{P}\{x - \varepsilon < \sqrt{n}\sigma_n Z/c_n < x+\varepsilon\} \ge \mathbb{P}\{x - \varepsilon < \sqrt{n}\sigma_n Z/c_n\}.$$

It is now evident that (3.14) implies (3.13) and the proof of the lemma is complete if $x \ge 0$. If $x < 0$, the lemma follows by symmetry. □

We are now ready to prove (2.13). By monotonicity of the exponential function and the definition of $\alpha_0$ we have

$$\sum_{n=1}^{\infty} n^{-1} \exp\left(-\frac{\alpha^2 c_n^2}{2n\sigma_n^2}\right) \begin{cases} = \infty, & \text{if } \alpha < \alpha_0, \\ < \infty, & \text{if } \alpha > \alpha_0. \end{cases}$$

Therefore if $\alpha_0 = \infty$, it trivially follows from Lemma 4 that $A \supset \mathbb{R}$, which of course implies that $A = [-\infty, \infty]$.

Assume now that $0 < \alpha_0 < \infty$. If $|x| \le \alpha_0$ and consequently $(|x| - \varepsilon)_+ < \alpha_0$, $\forall \varepsilon > 0$, we see that the series in Lemma 4 diverge for any $\varepsilon > 0$ so that $[-\alpha_0, \alpha_0] \subset A$.

Likewise, it follows that these series converge if $|x| > \alpha_0 \ge 0$ and $\varepsilon$ is sufficiently small. Thus such points are outside $A$ which implies that $A = [-\alpha_0, \alpha_0]$ and the first part of Theorem 3 has been proven.

If $\alpha_0 = \infty$, then (2.13) immediately implies (2.14), but if $\alpha_0$ is finite the lim sup in (2.14) still could be infinite. For that reason we have to add an extra argument to rule this out. Of course, we could apply Corollary 10 of [20], but since we already know the cluster set we do not need a precise upper bound for the lim sup. Once we know that the lim sup is finite it follows from (2.13) that it must be equal to $\alpha_0$. Here is a simple direct argument establishing this missing part of (2.14).

Proof of (2.14). We assume that $\alpha_0 < \infty$. Choosing $\delta = 1$, it follows that there exists an $\alpha > \alpha_0$ such that

$$\sum_{n=1}^{\infty} n^{-1} \exp(-\alpha^2 c_n^2/(2n\sigma_n^2)) < \infty, \tag{3.15}$$

where $\sigma_n^2 = H(c_n)$.

Set $n_k = 2^k$, $k \ge 1$. Then (3.15) immediately implies that

$$\begin{aligned}\sum_{n=1}^{\infty} n^{-1} \exp(-\alpha^2 c_n^2/(2n\sigma_n^2)) &\ge \sum_{k=1}^{\infty} \sum_{n=n_k}^{n_{k+1}-1} n^{-1} \exp(-\alpha^2 c_n^2/(2n\sigma_n^2)) \\ &\ge \sum_{k=1}^{\infty} \log(2) \exp(-\alpha^2 c_{n_{k+1}}^2/(2n_k\sigma_{n_k}^2)).\end{aligned}$$

Recalling (2.11) which implies that for some constant $K \ge 1$, $c_{n_{k+1}}/c_{n_{k-1}} \le 4K$, we find that

$$\sum_{k=2}^{\infty} \exp(-8K^2\alpha^2 c_{n_{k-1}}^2/(n_k\sigma_{n_k}^2)) < \infty. \tag{3.16}$$

We next employ Theorem 3 on page 74 in [3]. Assuming that the underlying probability space is rich enough and using (3.3), we can define a sequence of independent normal mean-zero random variables $Y_n$, $n \ge 1$, where $\mathrm{Var}(Y_n) =: \tilde{\sigma}_n^2 = \mathrm{Var}(XI\{|X| \le c_n\})$ so that we have for the sums $T_n = \sum_{i=1}^n Y_i$, $n \ge 1$, with probability 1,

$$(S_n - T_n)/c_n \to 0 \qquad \text{as } n \to \infty. \tag{3.17}$$

It is thus sufficient to prove that with probability 1,

$$\limsup_{n\to\infty} |T_n|/c_n \le 4K\alpha. \tag{3.18}$$

By the Borel–Cantelli lemma this follows once we have shown that

$$\sum_{k=1}^{\infty} \mathbb{P}\Big\{\max_{1\le m\le n_{k+1}} |T_m| \ge 4K\alpha c_{n_k}\Big\} < \infty. \tag{3.19}$$

But using a standard maximal inequality for normal random variables along with the fact that $\tilde{\sigma}_m^2 \le H(c_{n_{k+1}}) = \sigma_{n_{k+1}}^2, 1 \le m \le n_{k+1}$, we find that

$$\mathbb{P}\Big\{\max_{1\le m\le n_{k+1}} |T_m| \ge 4K\alpha c_{n_k}\Big\} \le 2\exp(-8K^2\alpha^2 c_{n_k}^2/(n_{k+1}\sigma_{n_{k+1}}^2)), \tag{3.20}$$

which in view of (3.16) implies (3.19). This completes the proof of (2.14). □

Proof of Theorem 4. Let $0 < d_n \le c_n$, $\sigma_{n,1}^2 := H(d_n)$ and define

$$\alpha_1 = \sup\Big\{\alpha \ge 0 : \sum_{n=1}^{\infty} n^{-1} \exp\Big(-\frac{\alpha^2 c_n^2}{2n\sigma_{n,1}^2}\Big) = \infty\Big\}.$$

As we have

$$\exp\Big(-\frac{\alpha^2 c_n^2}{2n\sigma_{n,1}^2}\Big) \le \exp\Big(-\frac{\alpha^2 c_n^2}{2n\sigma_n^2}\Big),$$

it is trivial that $\alpha_1 \le \alpha_0$.

We now consider normalizing sequences $c_n$ satisfying $a := \liminf_{n\to\infty} c_n/\gamma_n > 0$ and we choose $d_n \le c_n$ so that condition (2.15) is satisfied or, equivalently, $d_n = c_n/(Ln)^{\varepsilon_n}$, where $\varepsilon_n \to 0$. Let further $\Delta_n = \sigma_n^2 - \sigma_{n,1}^2$. In order to show that $\alpha_0 = \alpha_1$ it is enough to prove that

$$\sum_{n=1}^{\infty} n^{-1} \exp\Big(-\frac{\varepsilon c_n^2}{n\Delta_n}\Big) < \infty \qquad \forall \varepsilon > 0. \tag{3.21}$$

To see that, choose a $\delta > 0$, and observe that

$$\exp\left(-\frac{\alpha^2 c_n^2}{2n\sigma_n^2}\right) \le \exp\left(-\frac{\alpha^2 c_n^2}{2n(1+\delta)\sigma_{n,1}^2}\right) + \exp\left(-\frac{\alpha^2 c_n^2}{2n(1+\delta^{-1})\Delta_n}\right).$$

From (3.21) it is then obvious that $\alpha_0 \le \sqrt{1+\delta}\,\alpha_1$. Since we can choose $\delta$ arbitrarily small, we see that $\alpha_1 = \alpha_0$.

To prove that $\alpha_0 \le 1/a$, we set $d_n = c_n/(2aLLn)$ and use the fact that for large $n$ $\sigma_{n,1}^2 \le H(K(n/LLn)) \le K^2(n/LLn)LLn/n$. Replacing $\sigma_{n,1}^2$ in the definition of $\alpha_1$ by this upper bound, we readily obtain that $\alpha_0 = \alpha_1 \le 1/a$ and Theorem 4 has been proven subject to the verification of (3.21). □

Proof of (3.21). We use the same idea as in the proof of (2.14) of [5]. Recall that we have by definition of the $K$-function

$$H(K(x)) = \mathbb{E}X^2 I\{|X| \le K(x)\} \le K^2(x)/x, \qquad x > 0, \tag{3.22}$$

and

$$M(K(x)) = \mathbb{E}|X| I\{|X| \ge K(x)\} \le K(x)/x, \qquad x > 0. \tag{3.23}$$

To estimate $\Delta_n = \mathbb{E}X^2 I\{d_n < |X| \le c_n\}$ we observe that by (3.22) and Cauchy–Schwarz,

$$\begin{aligned}\Delta_n &\le \mathbb{E}X^2 I\{|X| \le K(n/LLn)\} + \mathbb{E}X^2 I\{|X| > K(n/LLn), |X| \le c_n\}\\ &\le K^2(n/LLn)LLn/n\\ &\quad + (\mathbb{E}|X| I\{|X| > K(n/LLn)\})^{1/2}(\mathbb{E}|X|^3 I\{|X| \le c_n\})^{1/2}.\end{aligned}$$

By assumption there exists an $n_0 \ge 1$ so that $c_n \ge aK(n/LLn)LLn$ for $n \ge n_0$. This implies in conjunction with (3.23),

$$\Delta_n \le \frac{c_n^2}{a^2 n}[(LLn)^{-1} + (na^3 \mathbb{E}|X|^3 I\{|X| \le c_n\}/c_n^3)^{1/2}], \qquad n \ge n_0. \tag{3.24}$$

Set $N_0 = \{n \ge n_0 : na^3 \mathbb{E}|X|^3 I\{|X| \le c_n\}/c_n^3 \le (LLn)^{-2}\}$. Then we have

$$\Delta_n \le 2c_n^2/(a^2 nLLn), \qquad n \in N_0. \tag{3.25}$$

As $d_n = c_n/(Ln)^{\varepsilon_n}$, where $\varepsilon_n \to 0$, we trivially have for $n \ge 1$,

$$\Delta_n \le \mathbb{E}|X|^3 I\{|X| \le c_n\}/d_n = (Ln)^{\varepsilon_n} \mathbb{E}|X|^3 I\{|X| \le c_n\}/c_n. \tag{3.26}$$

Employing the two bounds for $\Delta_n$ and recalling (3.1) we obtain via the trivial inequality $\exp(-x) \le 2x^{-1}\exp(-x/2)$ that

$$\begin{aligned}&\sum_{n \in N_0} n^{-1} \exp\{-\varepsilon c_n^2/(n\Delta_n)\}\\ &\quad \le 2\varepsilon^{-1} \sum_{n=1}^{\infty} (\Delta_n/c_n^2)(Ln)^{-\varepsilon a^2/4}\\ &\quad \le 2\varepsilon^{-1} \sum_{n=1}^{\infty} \mathbb{E}|X|^3 I\{|X| \le c_n\} c_n^{-3} (Ln)^{\varepsilon_n - \varepsilon a^2/4} < \infty.\end{aligned} \tag{3.27}$$

If $n \in N_1 = \{n \ge n_0 : n \notin N_0\}$, then $\Delta_n \le 2c_n^2\{\mathbb{E}|X|^3 I\{|X| \le c_n\}/(anc_n^3)\}^{1/2}$, and it follows from $e^{-x} \le 2/x^2$ that

$$\begin{aligned}&\sum_{n \in N_1} n^{-1} \exp\{-\varepsilon c_n^2/(n\Delta_n)\}\\ &\quad \le \sum_{n \in N_1} n^{-1} \exp(-(\varepsilon/2)\sqrt{a/n}(\mathbb{E}|X|^3 I\{|X| \le c_n\}/c_n^3)^{-1/2})\\ &\quad \le 8a^{-1}\varepsilon^{-2} \sum_{n=1}^{\infty} \mathbb{E}|X|^3 I\{|X| \le c_n\}/c_n^3 < \infty.\end{aligned} \tag{3.28}$$

This shows that (3.21) holds. □

Note that in the above proof we only use property (2.10) so that this relation holds for any sequence $c_n$ of positive real numbers such that $c_n/\sqrt{n}$ is nondecreasing. To conclude this section we show how the lower bound part of (1.8) follows from Theorem 3. To that end, it is sufficient to prove:

LEMMA 5. *If* $\limsup_{n\to\infty} c_n/\gamma_n \le b < \infty$, *we have* $\limsup_{n\to\infty} |S_n|/c_n \ge 1/b$ *a.s.*

PROOF. We apply Theorem 3 with $\delta = 1$. Then we have for large $n$,

$$\sigma_n^2 = H(c_n) \ge H(K(n/LLn)) + K(n/LLn)\mathbb{E}|X|I\{K(n/LLn) < |X| \le c_n\}$$

which by definition of the $K$-function is equal to

$$K^2(n/LLn)LLn/n - K(n/LLn)\mathbb{E}|X|I\{|X| > c_n\}.$$

Recalling (3.4) we see that $\liminf_{n\to\infty} n\sigma_n^2/\{K^2(n/LLn)LLn\} \ge 1$ which in turn implies that $\alpha_0 \ge 1/b$. □

**4. Proofs of Theorems 1 and 2.** We first note that by regular variation of $\Psi^{-1}$ we have

$$\limsup_{x\to\infty} \frac{\Psi^{-1}(xLLx)}{x^2LLx} H(x) = \limsup_{n\to\infty} \frac{nLLn}{a_n^2} H(a_n/LLn). \tag{4.1}$$

If one has a lower bound for the above lim sup one can infer that, along some subsequence, $\sigma_n^2 = H(a_n) \ge ch(n)/LLn$ for a positive $c$ which will imply that the series in the definition of $\alpha_0$ diverge for small positive $\alpha$ provided that the function $h$ is of very slow variation. This way we can prove that $\alpha_0$ is positive. (See Part 3 of the proof.)

If one has an upper bound for the above lim sup, one can in principle use the same approach to obtain an upper bound for $\alpha_0$. The problem here is that the above condition is not at the "natural" truncation level $a_n$. To overcome this difficulty, we first show (see Part 1) that under the assumptions of Theorem 1 we have $\liminf_{n\to\infty} a_n/\gamma_n > 0$ so that we can apply Theorem 4 which allows us to choose various truncation levels. Once this has been done, the upper bound (see Part 2) is straightforward (since any upper bound on a lim sup holds eventually in $n$).

It is then clear that the cluster set $A = C(\{S_n/a_n : n \ge 1\})$ is a bounded symmetric interval $[-\alpha_0, \alpha_0]$, and we shall show that $(1-q)^{1/2}\lambda \le \alpha_0 \le \lambda$, which clearly implies (2.2). As a matter of fact we then obtain a slightly stronger result, namely that under assumption (2.1) we have

$$(1-q)^{1/2}\lambda \le -\liminf_{n\to\infty} S_n/a_n = \limsup_{n\to\infty} S_n/a_n \le \lambda \qquad \text{a.s.}$$

Then Theorem 2 (with the extra information about the cluster set) is obvious and it is thus enough to prove Theorem 1.

*Part* 1. We need the following lemmas.

LEMMA 6. *Let $\Psi$ be as in Theorem 1. Assume that for some $\lambda \geq 0$,*

$$\limsup_{t\to\infty} \frac{\Psi^{-1}(tLLt)}{t^2 LLt} H(t) \leq \frac{\lambda^2}{2}. \tag{4.2}$$

*Then we also have*

$$\limsup_{t\to\infty} \frac{\Psi^{-1}(tLLt)}{tLLt} M(t) \leq C(1+\sqrt{2})\lambda^2, \tag{4.3}$$

*where $C > 0$ is a constant so that $\Psi^{-1}(x)/\Psi^{-1}(y) \leq C(x/y)^{3/2}$ for large $x \leq y$.*

PROOF. The existence of the constant $C$ follows easily from the Karamata representation of the slowly varying function $y \to \Psi^{-1}(y)/y^2$. (See, e.g., Theorem 1.3.1 in [1].) We thus can conclude that given $\delta > 0$, we have for large enough $t$

$$\begin{aligned}
M(t) &= \sum_{j=1}^{\infty} \mathbb{E}|X| I\{2^{j-1}t < |X| \leq 2^j t\} \leq \sum_{j=1}^{\infty} H(2^j t)/(2^{j-1}t) \\
&\leq (\lambda^2/2 + \delta) \sum_{j=1}^{\infty} LL(2^j t) 2^{j+1} t / \Psi^{-1}(2^j t LL(2^j t)) \\
&= (\lambda^2 + 2\delta) tLLt/\Psi^{-1}(tLLt) \\
&\quad \times \sum_{j=1}^{\infty} 2^j (LL(2^j t)/LLt)(\Psi^{-1}(tLLt)/\Psi^{-1}(2^j tLL(2^j t))) \\
&\leq C(\lambda^2 + 2\delta) tLLt/\Psi^{-1}(tLLt) \sum_{j=1}^{\infty} 2^j \{LLt/(LL(2^j t)\}^{1/2} 2^{-3j/2} \\
&\leq C(\lambda^2 + 2\delta) tLLt/\Psi^{-1}(tLLt) \sum_{j=1}^{\infty} 2^{-j/2} \\
&= C(1+\sqrt{2})(\lambda^2 + 2\delta) tLLt/\Psi^{-1}(tLLt).
\end{aligned}$$

Since $\delta$ can be chosen arbitrarily small, we obtain assertion (4.3). □

LEMMA 7. *Let $\Psi$ be as in Theorem 1. Then assumption (4.2) for some $\lambda \geq 0$ implies that $\liminf_{n\to\infty} \Psi(n)/K(n/LLn)LLn > 0$.*

PROOF. Recall that $a_n = \Psi(n)$. From Lemma 6 and assumption (4.2) it follows that there exists a positive constant $C'$ so that

$$\limsup_{t\to\infty} (H(t) + tM(t))\Psi^{-1}(tLLt)/(t^2 LLt) \leq C'\lambda^2 < \infty \tag{4.4}$$

which implies that for large enough $n$

$$G(a_n/LLn) \geq cn/LLn \tag{4.5}$$

where $1 \geq c > 0$, and, consequently,

$$a_n/LLn \geq K(cn/LLn) \geq cK(n/LLn) \tag{4.6}$$

and the lemma has been proven. □

REMARK 5. By a refinement of the above argument (where one has to choose the constant $c$ depending on $\lambda$ and show that $c$ goes to infinity as $\lambda$ goes to zero) one can also prove that if $\lambda = 0$ we have $\Psi(n)/K(n/LLn)LLn \to \infty$ as $n \to \infty$. Using this observation, one can infer the sufficiency part of Corollary 3 from (1.10).

*Part* 2 (*the upper bound*). In Section 2 we already have noted that the sequence $a_n$ satisfies assumption (2.10) and (2.11). Using the trivial fact that $\mathbb{E}\Psi^{-1}(|X|) < \infty$ if and only if $\sum_{n=1}^{\infty} \mathbb{P}\{|X| \geq a_n\} < \infty$, we see that Theorem 3 applies so that $\limsup_{n\to\infty} |S_n|/a_n = \alpha_0$ a.s. It remains to be shown that

$$\alpha_0 \leq \lambda. \tag{4.7}$$

In view of Part 1 we can apply Theorem 4 and it is sufficient to prove that if $\sigma_n^2 = H(a_n/LLn)$, we have

$$\sum_{n=1}^{\infty} n^{-1} \exp(-\alpha^2 h(n)/(2\sigma_n^2)) < \infty \qquad \forall\, \alpha > \lambda. \tag{4.8}$$

On account of (4.1) and (2.1), it follows that $2\sigma_n^2 \leq (\alpha-\delta)^2 h(n)/LLn$ for large $n$, where $\delta = (\alpha - \lambda)/2$.

This in turn implies $\exp(-\alpha^2 h(n)/(2\sigma_n^2)) \leq (Ln)^{-\eta^2}$, where $\eta = \alpha/(\alpha-\delta) > 1$. This clearly proves (4.8) and consequently (4.7).

*Part* 3 (*the lower bound and the converse to Theorem* 1). We present our last lemma from which we can infer both the lower bound in (2.2) and the converse to (2.3).

LEMMA 8. *Let $X : \Omega \to \mathbb{R}$ be a random variable satisfying for some $0 < \lambda < \infty$,*

$$\limsup_{x\to\infty} \frac{\Psi^{-1}(xLLx)}{x^2 LLx} H(x) \geq \lambda^2/2. \tag{4.9}$$

*If $h \in \mathcal{H}_q$ where $0 \leq q < 1$, we have*

$$\limsup_{n\to\infty} |S_n|/a_n \geq (1-q)^{1/2}\lambda \qquad a.s. \tag{4.10}$$

PROOF. It is sufficient to prove the lemma under the additional assumption

$$\mathbb{E}\Psi^{-1}(|X|) < \infty \quad \text{and} \quad \mathbb{E}X = 0. \tag{4.11}$$

To see that note that $\limsup_{n\to\infty} |S_n|/a_n < \infty$ a.s. implies that $\limsup_{n\to\infty} |X_n|/a_n < \infty$ a.s. By Kolmogorov's 0–1 law and the Borel–Cantelli lemma it then follows that

$$\sum_{n=1}^{\infty} \mathbb{P}\{|X| \geq a_n\} < \infty \tag{4.12}$$

which is equivalent to $\mathbb{E}\Psi^{-1}(|X|) < \infty$. So if this expectation is infinite, then by contraposition the lim sup in (4.10) is infinite.

By the strong law of large numbers this is also the case if $\mathbb{E}X \neq 0$.

Finally, without loss of generality, we can assume that $a_n/\sqrt{n} \nearrow \infty$. [Note that Theorem 3 with $c_n = n^{1/2}(LLn)^{1/3}$ trivially implies that $\limsup_{n\to\infty} |S_n|/a_n = \infty$ a.s. if $\mathbb{E}X = 0$ and (4.12) is satisfied with $a_n = O(\sqrt{n})$.]

Under the above assumptions Theorem 3 applies. We shall show that

$$\alpha_0 \geq (1-q)^{1/2}\lambda.$$

It then follows that $[-(1-q)^{1/2}\lambda, (1-q)^{1/2}\lambda] \subset A = C(\{S_n/a_n; n \geq 1\})$. This trivially implies (4.10). By definition of $\alpha_0$ and monotonicity, it is enough to prove that

$$\sum_{n=1}^{\infty} n^{-1} \exp(-\alpha^2 h(n)/(2\sigma_n^2)) = \infty, \qquad 0 < \alpha < (1-q)^{1/2}\lambda, \tag{4.13}$$

where $\sigma_n^2 = H(a_n/LLn)$.

Recalling (4.1) we see that

$$\limsup_{n\to\infty} \{LLn/h(n)\}\sigma_n^2 \geq \lambda^2/2. \tag{4.14}$$

Given an $\alpha$ as above, choose $0 < \tau' < 1-q$ so that $\alpha^2 = \tau'\lambda^2$ and set $\tau = \tau' + \delta/2$, where $\delta = 1 - q - \tau'$. Let $f_\tau$ be defined as in Section 2. On account of (4.14) and the definition of $\mathcal{H}_q$ we can find a subsequence $m_k \nearrow \infty$ so that

$$\sigma_{m_k}^2 \geq \frac{\lambda^2}{2}\left(1 - \frac{1}{k}\right)\frac{h(m_k)}{LLm_k} \tag{4.15}$$

and

$$h(m_k) \geq (1 - 1/k)h(m_k f_\tau(m_k)), \qquad k \geq 1. \tag{4.16}$$

Combining the last two relations we readily obtain by monotonicity of $\sigma_n^2$ in $n$ that

$$\sigma_n^2 \geq \frac{\lambda^2}{2}\left(1 - \frac{1}{k}\right)^2 \frac{h(n)}{LLn}, \qquad m_k \leq n \leq n_k := [m_k f_\tau(m_k)], \tag{4.17}$$

which in turn implies that

$$\sum_{n=m_k}^{n_k} n^{-1}\exp\Big(-\frac{\alpha^2 h(n)}{2\sigma_n^2}\Big) \geq \log\Big(\frac{n_k+1}{m_k}\Big)(Ln_k)^{-\tau'/(1-1/k)^2}. \tag{4.18}$$

As we have $\log f_\tau(m_k) = (\log m_k)^\tau \leq \log m_k$ we get for large $k$

$$\log\Big(\frac{n_k+1}{m_k}\Big)(Ln_k)^{-\tau'/(1-1/k)^2} \geq (Lm_k)^\tau (2Lm_k)^{-(\tau'+\delta/4)} \geq (Lm_k)^{\delta/4}/2$$

which goes to infinity. Recalling (4.18), we obtain (4.13) and the lemma has been proven. □

Combining (4.7) and Lemma 8, we obtain (2.2). Moreover, in the proof of Lemma 8 we have already shown that the assumptions $\mathbb{E}\Psi^{-1}(|X|) < \infty$ and $\mathbb{E}X = 0$ are necessary for (2.2) to hold.

Furthermore, if $\limsup_{x\to\infty} \frac{\Psi^{-1}(xLLx)}{x^2LLx} H(x) = \infty$ and if $q < 1$ we can infer from Lemma 8 (with arbitrarily large $\lambda$) that

$$\limsup_{n\to\infty} |S_n|/a_n = \infty.$$

This clearly shows that (2.1) for some $\lambda < \infty$ is necessary for (2.3) to hold.

Likewise, if $\limsup_{x\to\infty} \frac{\Psi^{-1}(xLLx)}{x^2LLx} H(x) = 0$, we obtain by (4.7) (with $\lambda = 0$) $\limsup_{n\to\infty} |S_n|/a_n = 0$ a.s. and it is clear that this lim sup can only be positive if condition (2.1) holds for some $\lambda > 0$.

**5. Further examples.** We first give a corollary to Theorem 1 where $h \in \mathcal{H}_q$ and $0 < q < 1$. We consider

$$h_q(x) = \exp\{(Lx)^q\} \quad \text{and} \quad \Psi_q(x) = \sqrt{x\exp\{(Lx)^q\}}.$$

It is easy to see that $h_q \in \mathcal{H}_q$. Write

$$H_q(x) = \frac{x^2}{\exp\{2^q(Lx)^q\}}.$$

One can check that

$$\lim_{x\to\infty} \frac{\Psi_q(H_q(x))}{x} = \begin{cases} 1, & \text{if } 0 < q < 1/2,\\ e^{-1/4}, & \text{if } q = 1/2.\end{cases}$$

We thus have

$$\lim_{x\to\infty} \frac{\Psi_q^{-1}(x)}{H_q(x)} = \begin{cases} 1, & \text{if } 0 < q < 1/2,\\ e^{1/2}, & \text{if } q = 1/2.\end{cases}$$

For $1/2 < q < 1$, the precise asymptotic expansion of $\Psi_q^{-1}(x)$ is a little bit complicated and is left to the interested reader. Applying Theorem 1 to the case where $0 < q \leq 1/2$, we have the following result.

Corollary 4. *Let $0 < q \le 1/2$. If there exists a constant $0 \le \lambda < \infty$ such that*

$$\mathbb{E}X = 0, \qquad \mathbb{E}\left(\frac{X^2}{\exp(2^q(L|X|)^q)}\right) < \infty \tag{5.1}$$

*and*

$$\begin{aligned} \limsup_{x\to\infty} \frac{LLx}{\exp(2^q(Lx)^q)} H(x) &= \frac{\lambda^2}{2} \qquad \text{if } 0 < q < 1/2, \\ \limsup_{x\to\infty} \frac{e^{1/2}LLx}{\exp(\sqrt{2Lx})} H(x) &= \frac{\lambda^2}{2} \qquad \text{if } q = 1/2, \end{aligned} \tag{5.2}$$

*then we have with probability* 1,

$$(1-q)^{1/2}\lambda \le \limsup_{n\to\infty} \frac{|S_n|}{\sqrt{n\exp((Ln)^q)}} \le \lambda.$$

Here of course it would be interesting to know whether our bounds for the above lim sup are sharp. In principle one can calculate the precise value of the lim sup via Theorem 3 and it may depend on the distribution of $X$. One might wonder whether all values in the interval $[(1-q)^{1/2}\lambda, \lambda]$ can occur or whether one can improve the general lower bound we have found.

Let us take another look at Theorem 1. For a given sequence of i.i.d. mean-zero random variables $\{X, X_n; n \ge 1\}$, we may want to know if there exists a sequence of positive real numbers $\{a_n; n \ge 1\}$ such that

$$0 < \limsup_{n\to\infty} \frac{|S_n|}{a_n} < \infty \qquad \text{a.s.} \tag{5.3}$$

holds and if it does, how to find it. To answer this question, we may try the following method. Let again $H(x) = \mathbb{E}(X^2 I\{|X| \le x\})$, $x \ge 0$, and suppose there exists a positive and nondecreasing slowly varying function $\phi(x)$ such that

$$\limsup_{x\to\infty} \frac{H(x)}{\phi(x)} = 1; \tag{5.4}$$

we then take $\Psi(x)$ such that

$$\frac{\Psi^{-1}(xLLx)}{x^2 LLx}\phi(x) = 1. \tag{5.5}$$

Thus, $\Psi(x)$ satisfies

$$\Psi^{-1}(x) \sim \frac{x^2}{(LLx)\phi(x/LLx)} \qquad \text{as } x \to \infty, \tag{5.6}$$

which is equivalent to

$$\Psi(x) \sim (x\phi(\Psi(x)/LLx)LLx)^{1/2} \qquad \text{as } x \to \infty. \tag{5.7}$$

If $h(x) = \phi(\Psi(x)/LLx)LLx \in \mathcal{H}_q$, where $q < 1$, then (5.3) holds with $a_n = \Psi(n)$ if and only if

$$\mathbb{E}\left(\frac{X^2}{\phi(|X|/LL|X|)LL|X|}\right) < \infty. \tag{5.8}$$

Of course, if $H$ is already slowly varying at infinity, which implies that $X$ is in the domain of attraction to the standard normal distribution, then this result holds in general if we choose $H = \phi$ without assuming that $q < 1$. This follows, for instance, from Theorem 1 of [16]. But even in this situation it can be very helpful to work with a different (and larger) slowly varying function $\phi$. To demonstrate this we shall look at an example which was also discussed by Feller [6] and Pruitt ([22], Example 9.4).

Example. Let $\{X, X_n; n \ge 1\}$ be a sequence of real-valued i.i.d. random variables with the common symmetric probability density function

$$f(x) = \frac{1}{|x|^3} I\{|x| \ge 1\}.$$

For this example, Pruitt ([22], page 44) pointed out that it would be possible to find a normalizing sequence $\{a_n; n \ge 1\}$ such that

$$\limsup_{n\to\infty} \frac{|S_n|}{a_n} = 1 \qquad \text{a.s.}$$

Can the normalizing sequence $\{a_n; n \ge 1\}$ be explicitly given? Pruitt [22] did not answer this question but mentioned that it would not be a very nice normalizing sequence. Using our procedure above, we can find a normalizing sequence of the form $\sqrt{nh(n)}$ with $h$ slowly varying which is not as unreasonable as one might expect. In fact, for this example, $H(x) = 2Lx$, $x \ge 0$. If $\phi_1(x) = 2Lx$, $x \ge 0$ is chosen to be the $\phi(x)$, then by (5.7),

$$\Psi_1(x) \sim (xLxLLx)^{1/2} \qquad \text{as } x \to \infty.$$

It is easy to check that (5.8) does not hold with $\phi_1(x) = 2Lx$ which implies

$$\limsup_{n\to\infty} \frac{|S_n|}{(nLnLLn)^{1/2}} = \infty \qquad \text{a.s.}$$

However, we may choose $\phi_2(x) = 2Lx(1 + LLx\sin^2(LLLx))$, $x \ge 0$, to be the $\phi(x)$. It is easily checked that $\phi_2'(x) \ge 0$ so that this is a function in $\mathcal{H}$. After some calculation it also follows that $\phi_2 \in \mathcal{H}_0$. We obviously have $\limsup_{x\to\infty} H(x)/\phi_2(x) = 1$. Moreover, using that $\phi_2(x/LLx) \sim \phi_2(x)$ as $x \to \infty$, we infer from (5.7)

$$\Psi_2(x) \sim \left(\frac{x}{2}\phi_2(x)LLx\right)^{1/2} \qquad \text{as } x \to \infty$$

and (5.8) holds with $\phi_2(x)$ since

$$\begin{aligned}\mathbb{E}\left(\frac{X^2}{\phi_2(|X|)LL|X|}\right) &= 2\int_1^C \frac{1}{x\phi_2(x)LLx}\,dx + 2\int_C^\infty \frac{1}{x\phi_2(x)LLx}\,dx\\ &= 2\int_1^C \frac{1}{x\phi_2(x)LLx}\,dx + \int_1^\infty \frac{1}{1+e^y\sin^2 y}\,dy\\ &< \infty,\end{aligned}$$

where $C = e^{e^e}$. Thus, by Theorem 2, we have

$$\limsup_{n\to\infty} \pm S_n/\sqrt{2nLnLLn(1+LLn\sin^2(LLLn))} = 1 \qquad \text{a.s.}$$

**Acknowledgments.** The authors thank the referee for very helpful suggestions to improve the presentation and for mentioning the work of Martikainen. They are grateful to Professors R. J. Tomkins and D. M. Mason for reading the manuscript and for many valuable comments.

Departement Wiskunde
Vrije Universiteit Brussel
Pleinlaan 2
B-1050 Brussel
Belgium
e-mail: ueinmahl@vub.ac.be

Department of Mathematical Sciences
Lakehead University
Thunder Bay, Ontario
Canada P7B 5E1
e-mail: lideli@giant.lakeheadu.ca